\documentclass{article}
\usepackage{amssymb,latexsym,amsmath,epsfig,amsthm,hyperref,tipa,comment}

\newtheorem*{main}{Theorem}

\newtheorem{howmanysmoothers}{Lemma}
\newtheorem{ssize}[howmanysmoothers]{Lemma}
\newtheorem{maximums}[howmanysmoothers]{Lemma}

\begin{document}

\vspace*{-2cm}

\Large
 \begin{center}
Smooth sums with small spacings \\ 

\hspace{10pt}

\large
Wouter van Doorn and Anneroos R.F. Everts \\

\hspace{10pt}

\end{center}

\hspace{10pt}

\normalsize

\vspace{-25pt}
\centerline{\bf Abstract}
Solving a problem by Erd\H{o}s, we prove that every positive integer $n$ can be written as a sum $$n = b_{1} + b_{2} + \ldots + b_{r}$$ of distinct $3$-smooth integers with $1 \le b_{1} < b_{2} < \ldots < b_{r} < 6b_{1}$.

\section{Introduction}
Let $A = (a_1, a_2, \ldots)$ be the infinite increasing sequence of $3$-smooth integers. That is, for every index $i$ there are non-negative integers $x_i, y_i$ for which $a_i = 2^{x_i} 3^{y_i}$, while $a_{i+1} > a_i$ for all $i \in \mathbb{N}$. In the early $1990$s, Erd\H{o}s conjectured that every positive integer $n$ can be written as a sum of distinct $a_i$ such that no summand divides another. As was quickly realized however (even before it was written down anywhere in the literature, for the first time in \cite{Er92b}), this conjecture actually has a very short induction proof. Indeed, one can make the stronger induction hypothesis that, for all even $n$, all summands are even as well. For even $n$ we are then done by applying the hypothesis to $\frac{n}{2}$, while for odd $n$ one can apply the induction hypothesis to $n - 3^{\lfloor \log_3 n \rfloor}$. \\

In general, sequences such that every large enough integer can be written as a sum of distinct elements where no summand divides another, are called $d$-complete sequences. They have been studied by Erd\H{o}s-Lewin \cite{ErLe96}, Ma-Chen \cite{MaCh16}, and Chen-Wu \cite{ChYu23}, mostly in the context of proving, for various fixed $a, b, c$, that the sequence of integers of the form $a^x b^y c^z$ is $d$-complete. \\

Going back to $3$-smooth integers, Blecksmith, McCallum and Selfridge show in \cite{BMS98} that $A$ remains $d$-complete even after removing finitely many elements from $A$, thereby proving the final conjecture from \cite{ErLe96}. Their proof actually provides a fully explicit procedure to write an integer $n$ as a sum of distinct and large $a_i$, and a variation of this procedure will feature prominently in this paper. \\

Note that no two summands $b_i, b_j$ can divide one another if $b_i < b_j < 2b_i$. Therefore, a stronger conjecture in this regard would be whether, for some $C \le 2$, every positive integer $n$ can be written as a sum $n = b_{1} + b_{2} + \ldots + b_{r}$ with $b_i \in A$ for all $i$, and $b_{1} < b_{2} < \ldots < b_{r} < Cb_1$. This stronger conjecture turns out to be false however (although it plausibly does hold if we consider $5$-smooth integers instead), as was already remarked in \cite{ErLe96}. A natural follow-up question is then if a constant $C > 2$ exists for which such a representation is possible for all $n \in \mathbb{N}$, even though we are by now leaving the realm of $d$-completeness. \\

The possible existence of such a $C$ was first considered in \cite{Er92b}, where Erd\H{o}s initially thought that `surely almost all integers cannot be written in this form'. This question was then repeated, much more neutrally, in \cite{ErLe96}, and it is now listed as Problem $845$ at Bloom's website \cite{bloom}. There, in the comment section, Cambie suggested that such a constant actually does exist, and checked with a computer that $C = \frac{32}{9} = 3.55\cdots$ works for all $n \le 10^5$. \\

In this paper we will more generally consider the sequence $A_p$ of positive integers of the form $2^x p^y$. We will then show that for all odd integers $p > 1$ there exists a constant $C_p$ such that every positive integer $n$ can be written as a sum of distinct elements of $A_p$ for which the ratio between any two summands is smaller than $C_p$. We do this by tweaking and generalizing the procedure from \cite{BMS98}, where for $p = 3$ we in particular obtain $C_p = 6$.

\section{Main result}
Let $p > 1$ be an odd integer, and define $A_p = (a_1, a_2, \ldots)$ as the infinite increasing sequence of all integers that can be written as $a_i = 2^{x_i} p^{y_i}$, for some non-negative integers $x_i$ and $y_i$. With $\log_2 x$ denoting the logarithm function to base $2$, we furthermore define the functions $f_0(x) = x$, $f_k(x) = \max\left(1, \left \lfloor \log_2 f_{k-1}(x)\right \rfloor\right)$ for $k \ge 1$, and $F(x) = \prod_{k \ge 0} f_k(x)$. Our main result can then be stated as follows.

\begin{main} \label{main}
For every odd integer $p > 1$ a constant $C_p$ exists such that every positive integer $n$ can be written as a sum $$n = b_{1} + b_{2} + \ldots + b_{r}$$ with $b_i \in A_p$ for all $i$ and $b_{1} < b_{2} < \ldots < b_{r} < C_pb_1$. \\

In general the constant $C_p$ can be taken to be equal to $\frac{1}{2}F(4p)$, and if either $p-1$ or $p+1$ is a power of two, then one can take $C_p = 2p$ or $C_p = 2(p+1)$ respectively. On the other hand, one cannot replace $C_p$ by any constant smaller than $p$.
\end{main}

\begin{proof}
We will first prove that $C_p < p$ would not be admissible. So let us assume that $C$ is a constant with $1 < C < p$, choose $\delta > 0$ and $\epsilon > 0$ sufficiently small so that $C(1+\epsilon) < p - \delta$, and let $N$ be large enough. For ease of reference, let us call a sum $b_{1} + b_{2} + \ldots + b_{r}$ of distinct elements from $A
_p$ with $b_{1} < b_{2} < \ldots < b_{r} < Cb_1$ \textit{short}. The goal is to show that the number of short sums with all elements smaller than or equal to $N$, is smaller than $N$. \\

For a non-negative integer $j$, let $x_j$ be equal to $(1+\epsilon)^j$, and define the interval $I_j = [x_j, x_{j+1})$. For any short sum with $b_{1} \in I_j$, we then get 
\begin{align*}
x_j &\le b_{1} \\
&\le b_{r} \\
&< Cb_{1} \\
&< Cx_{j+1} \\
&= C(1+\epsilon)x_j \\
&< (p-\delta)x_j. 
\end{align*}

In particular, with $X_j$ the number of elements of $A_p$ contained in the interval $[x_j, (p-\delta)x_j)$, the number of short sums with $b_1 \in I_j$ is at most $2^{X_j}$. From the discussion in Lecture $5$ from \cite{HaRa78}, we moreover have the following lemma bounding $X_j$.

\begin{howmanysmoothers} \label{howmany}
There exists a constant $c_p$ such that $X_j < \frac{\log x_j \log(p - \delta)}{\log 2 \log p} + c_p$ for all $j \ge 0$.
\end{howmanysmoothers}

Applying Lemma \ref{howmany} and using the fact that every subset sum with all elements smaller than or equal to $N$ must have $b_1 \in I_j$ for some $j \le L := \lfloor \frac{\log N}{\log (1+\epsilon)} \rfloor$, the total number of short sums at most $N$ is upper bounded by
\begin{align*}
\sum_{j=0}^L 2^{X_j} &< 2^{c_p} \sum_{j=0}^L x_j^{\frac{\log (p - \delta)}{\log p}} \\
&\le 2^{c_p} (L+1) x_L^{\frac{\log (p - \delta)}{\log p}} \\
&\le 2^{c_p} (L+1) N^{\frac{\log (p - \delta)}{\log p}}.
\end{align*}

Since this latter quantity is significantly smaller than $N$ when $N$ is large enough, we conclude that almost all positive integers cannot be represented as a short sum. \\

We may now focus on the other direction and prove that such representations do exist for all $n$, if we choose $C_p$ sufficiently large. Hence, from here on out we let $n$ be any arbitrary, but fixed, positive integer. \\

Let $S \subset A_p \setminus\{1\}$ be a finite set for which the set of subset sums of $S$ contains $|S|+1$ consecutive integers. To give an example, if $p = 19$, then we claim that we can take $S = \{2, 4, 16, 19\}$. One can verify that
\begin{align*}
18 &= 2 + 16, \\
19 &= 19, \\ 
20 &= 4 + 16, \\
21 &= 2 + 19 \text{ and} \\
22 &= 2 + 4 + 16,
\end{align*}

so that there are indeed $|S| + 1 = 5$ consecutive integers that can be written as a sum of distinct elements of $S$. \\

In general, from the conditions on $S$ it may not be immediately obvious what such sets look like or if they even exist. Luckily, their existence is not too hard to show.

\begin{ssize} \label{ssize}
For all odd integers $p > 1$ a set $S \subset A_p \setminus \{1\}$ exists with the following four properties:

\begin{enumerate}
	\item The set $S$ contains all even powers of two smaller than or equal to $|S|$.
	\item A positive integer $M_0 \le p$ exists such that all integers $x$ with $M_0 \le x \le M_0 + |S|$ can be written as a sum of distinct elements of $S$.
	\item The cardinality of $S$ is at most $1 + \lceil \log_2 p \rceil$.
	\item For the largest element $\max S$ of $S$ we have $p \le \max S \le 2^{\lceil \log_2 p \rceil}$.
\end{enumerate}
\end{ssize}

\begin{proof}
We claim that we can always take $S = \{p, 2, 4, \ldots, 2^{\lceil \log_2 p \rceil}\}$. For this set, the only non-trivial property is the second one. To prove that the second property is also satisfied, we note that 
\begin{align*}
p + |S| &= p + 1 + \lceil \log_2 p \rceil \\
&\le 2p \\
&< 2^{\lceil \log_2 p \rceil + 1}.
\end{align*}

This implies that all even integers smaller than or equal to $p + |S|$ can be written as a subset sum of $S \setminus\{p\}$ using their binary expansion, while all odd integers $x$ with $p \le x \le p + |S|$ can be written as $x = p + \sum_{y \in Y} y$ for some subset $Y \subseteq S \setminus\{p\}$ using the binary expansion of $x-p$. We conclude that the second property is satisfied with $M_0 = p$.
\end{proof}

So from now on, let $S$ and $M_0$ be such that the four conditions of Lemma \ref{ssize} all hold, and let $m$ be the smallest index with $a_1 + a_2 + \ldots + a_m > \frac{n}{M_0}$. In fact, we will need a couple of additional definitions. \\

We define $M_1 = M_0+|S|$, $M_2 = |S|$ and for $k \ge 3$, if $M_{k-1} > 1$, we further define $M_{k} = \lfloor \log_2 M_{k-1} \rfloor$. We note that only finitely many $M_k$ can exist, and we denote by $K$ the largest index for which $M_K$ is defined. The sequence $(u_1, u_2, \ldots, u_K)$ is then defined by $u_1 = 1$, $u_2 = \max S$ and $u_k = 2^{M_{k}}$ for $3 \le k \le K$. We now set $P_k = u_1u_2\cdots u_k$ for $1 \le k \le K$, and choose $C_p$ to be equal to the product $P_K$. Finally, we define the intervals $I_1 = [1, a_m)$ and $I_k = [P_{k-1}a_m, P_k a_m)$ for $2 \le k \le K$, and the sequence $(v_1, v_2, \ldots, v_{K})$ where $v_k$ is the index for which $a_{v_k+1} = P_{k} a_m$. That is, $a_{v_k}$ is the largest element of $A_p$ in the interval $I_k$. \\

Now, recalling the definition of $m$, we write 
\begin{equation} \label{binary}
n - M_0(a_1 + a_2 + \ldots + a_{m-1}) = a_{j_1} + a_{j_2} + \ldots + a_{j_{s}} 
\end{equation}

in binary, i.e. the $a_{j_i}$ on the right-hand side of equation (\ref{binary}) are distinct powers of two. \\

By adding the sums $M_0(a_1 + a_2 + \ldots + a_{m-1})$ and $a_{j_1} + a_{j_2} + \ldots + a_{j_{s}}$, we obtain the following representation of $n$:
\begin{equation} \label{firstrep}
n = c_1a_1 + c_2a_2 + \ldots + c_{k}a_{m-1} + c_{m}a_{m} + \ldots + c_{v_{K}}a_{v_{K}}.
\end{equation}

Here, $c_i \in \{M_0, M_0 + 1\}$ for $1 \le i \le v_1 = m-1$. For $i > v_1$ we have $c_i \in \{0, 1\}$ with $c_i = 1$ if, and only if, $a_i \ge a_m$ occurred as some power of two on the right-hand side of (\ref{binary}). By the definition of $m$ one can check that either side of (\ref{binary}) is smaller than 
\begin{align*}
M_0a_m &\le pa_m \\
&\le u_2a_m \\
&= P_2a_m,
\end{align*}

so that $c_i = 0$ for all $i > v_2$. \\

Using a variation on the procedure laid out as `the midgame' in \cite{BMS98}, we are going to transform the representation from equation (\ref{firstrep}) into a different representation of $n$ in such a way that, at the end, $c_i \in \{0, 1\}$ for all $i$, and where $c_i = 1$ implies $v_1 < i \le v_K$. Stated differently, $c_i$ will eventually be equal to $0$ unless $a_m \le a_i < C_pa_m$, which would finish the proof. \\

In step $i$ of the transformation procedure, we consider the coefficient $c_i$. If $c_i > 1$, we write either $c_i$ or $c_i-1$ as a sum of distinct elements of $A_p$. We then lower $c_i$ to either $0$ or $1$, while increasing $c_{i'}$ for certain $i' > i$, in such a way that equality in (\ref{firstrep}) is maintained. To elaborate on this, let us first assume $i \le v_1$. \\

In that case we know by (\ref{firstrep}) that we initially have $c_i \ge M_0$. If, moreover, $c_i \le M_1$, then we can write $c_i$ as a sum $a_{i,1} + a_{i,2} + \ldots + a_{i,t}$ of distinct elements of $S \subset A_p \setminus\{1\}$, by the second property of $S$ mentioned in Lemma \ref{ssize}. The term $c_ia_i$ in (\ref{firstrep}) can then be written as $a_{i,1}a_i + a_{i,2}a_i + \ldots + a_{i,t}a_i$. Since $A_p$ is multiplicatively closed, for every $j$ with $1 \le j \le t$ we have that $a_{i,j}a_i$ is equal to $a_{i'}$ for some $i' > i$. By decreasing $c_i$ to $0$ and increasing $c_{i'}$ by $1$ for all $i'$ for which $a_{i'}$ is equal to $a_{i,j}a_i$ for some $j$, equality in (\ref{firstrep}) is maintained. \\

We claim that $c_i \le M_1$ does indeed hold for $i \le v_1$, so that the above procedure works for all $i \le v_1$. To see this, consider which $i' < i$ can be responsible for increasing $c_i$. This can only happen if $\frac{a_i}{a_{i'}} \in S$, in which case it is possible that $c_i$ increased by $1$ in step $i'$. Hence, when we reach step $i$, $c_i$ has been increased by at most $|S|$ from its starting value. Now we remark that $c_i$ for $i \le v_1$ can have two different starting values; $M_0$ or $M_0 + 1$. In the first case, $c_i$ will be at most $M_0 + |S| = M_1$ when we reach step $i$. In the second case, we know that $a_i$ is a power of two that occurred in the binary expansion on the right-hand side of (\ref{binary}). With $a_i$ a power of two, $\frac{a_i}{a_{i'}}$ is of course even for all $a_{i'}$ dividing $a_i$. On the other hand, $S$ must contain at least one odd integer by the second property of Lemma \ref{ssize}, implying that $c_i$ has been increased by at most $|S|-1$ by the time we reach step $i$. Therefore, in this case we have that $c_i$ will be at most $M_0 + 1 + |S| - 1 = M_1$ as well. \\

Now let us explain how to transform $c_i$ if $i > v_1$. That is, when $a_i \notin I_1$. In this case, if $c_i \in \{0, 1\}$, then we do not do anything and we simply go to step $i+1$. On the other hand, if $c_i > 1$, then either $c_i$ or $c_{i}-1$ can be uniquely written as a sum $a_{i,1} + a_{i,2} + \ldots + a_{i,t'}$ of distinct even powers of two. Similarly to what we had in $I_1$, by decreasing $c_i$ to $c_i - (a_{i,1} + a_{i,2} + \ldots + a_{i,t'}) \in \{0, 1\}$ and increasing $c_{i'}$ by $1$ for all $i'$ for which $a_{i'}$ is equal to $a_{i,j}a_i$ for some $j$, we once again maintain equality in (\ref{firstrep}). And when at the end of step $i$ we have $c_{i'} \in \{0, 1\}$ for all $i' \ge i$, we stop. \\

Once we have finished step $v_k$, we see that all coefficients in $I_k$ have been brought down to either $0$ or $1$. This process certainly terminates at some point; in order for equality (\ref{firstrep}) to hold, we must have $c_i = 0$ at all times for all $i$ with $a_i > n$, so we are guaranteed to stop at or before step $n$. We claim however that this process already finishes in at most $v_{K-1}$ steps with $c_i$ still at $0$ for all $i > v_K$. This follows from the following lemma.

\begin{maximums} \label{maxx}
Throughout the entire procedure we have 

\begin{equation} \label{maxmk}
\max_{v_{k-1} < i \le v_k} c_i \le M_k 
\end{equation}

for all $k$ with $1 \le k \le K$. And if $c_i$ with $i > v_{k-1}$ was increased in step $i'$ for some $i' < i$, then $i' > v_{k-2}$.
\end{maximums}

\begin{proof}
Earlier we already proved (\ref{maxmk}) for $k = 1$, where the inequality $v_{k-1} < i$ here may be ignored. Following that same argument and applying the first property of $S$, we also get $c_i \le |S| = M_2$ for all $v_1 < i \le v_2$. Furthermore (where we again ignore inequalities with undefined terms), the second claim of Lemma \ref{maxx} is vacuous for $k \le 2$. Now we use induction on $i$, so assume that we are at the start of step $i$ with $v_{k-1} < i \le v_k$ for some $k \ge 3$. Further assume that $c_{i'} \le M_{k-1}$ for all $v_{k-2} < i' \le v_{k-1}$, $c_{i'} \le M_k < M_{k-1}$ for all $v_{k-1} < i' < i$, and that if $c_i$ was increased in step $i'$, then $i' > v_{k-2}$. In particular, if $c_i$ was increased in step $i'$, then $\frac{a_i}{a_{i'}}$ must be an even power of two smaller than or equal to $M_{k-1}$. From this we indeed conclude that $c_i$ is at most $\lfloor \log_2 M_{k-1} \rfloor = M_k$. \\

Moreover, when we then perform step $i$ and write either $c_i \le M_k$ or $c_i - 1 < M_k$ as a sum of even powers of two, then we can only increase $c_{i'}$ for some $i' > i$ if $\frac{a_{i'}}{a_{i}}$ is a power of two smaller than or equal to $M_k$, i.e. at most $2^{\lfloor \log_2 M_k \rfloor} = 2^{M_{k+1}} = u_{k+1}$. As the assumption $i \le v_k$ is equivalent to $a_i < P_ka_m$, we deduce $a_{i'} < u_{k+1}P_ka_m = P_{k+1}a_m$, and we see that in step $i$ we can only increase $c_{i'}$ if $i' \le v_{k+1}$.
\end{proof}

Finally, we need to prove that $C_p$ is smaller than or equal to the claimed values in the statement of the theorem. If $p+1$ is a power of two, we take $S = \{2, p, p+1\}$, which covers the integers $p$, $p+1$, $p+2$ and $p+3$ with its subset sums. This gives $M_2 = |S| = 3$ and $M_3 = 1$, so that with $u_2 = \max S = p+1$ and $u_3 = 2$ we obtain $C_p = 2(p+1)$. If $p-1$ is a power of two, we take $S = \{2, p-1, p\}$ instead, which covers the integers $p-1$, $p$, $p+1$ and $p+2$ with its subset sums. This still gives $M_2 = 3$ and $M_3 = 1$, while $u_2 = p$ and $u_3 = 2$ in this case, implying $C_p = 2p$. \\

Before we continue with the general case, we note that the set $S = \{2, p-1, p\}$ in the previous example is actually a multiset for $p = 3$. This implies that, when writing $c_i$ as a sum of distinct elements of $S$, we may use the element $2$ twice. This furthermore implies that, if we write $c_i = 4$ as $a_{i,1} + a_{i,2} = 2+2$ for example, then for $i'$ such that $a_{i'} = 2a_i$, we need to increase $c_{i'}$ twice. Even though this does not change anything in (our analysis of) the transformation procedure, it still seems worth it to be aware of this possibility, and we will come back to it shortly. \\

In the general case we have 
\begin{align*}
u_2 &= \max S \le 2^{\lceil \log_2 p \rceil} < 2p, \\
u_3 &= 2^{M_3} = 2^{\lfloor \log_2 M_2 \rfloor} \le M_2 = |S| \text{ and}\\
u_4 &= 2^{M_4} = 2^{\lfloor \log_2 M_{3} \rfloor} \le M_{3} = \lfloor \log_2 M_{2} \rfloor = f_1(|S|).
\end{align*}

The final equality generalizes via induction to $$u_k = 2^{M_k} \le \lfloor \log_2 M_{k-2} \rfloor = f_{k-3}(|S|)$$ for all $4 \le k \le K$. We therefore finally conclude 
\begin{align*}
C_p &= u_2u_3 \cdots u_K \\
&< 2p \prod_{k=0}^{K-3} f_k(|S|) \\
&< 2p \prod_{k=0}^{K-3} f_k(\lfloor \log_2 4p \rfloor) \\
&= 2p \prod_{k=1}^{K-2} f_k(4p) \\
&\le \frac{1}{2} F(4p). \qedhere
\end{align*}
\end{proof}

\section{A final optimization}
We recall that for $p = 3$ we needed $S = \{2, 2, 3\}$ to be a multiset, in order to end up with a better value of $C_3$. Indeed, one can check that $S = \{2, 3, 4\}$ works as well, but would have given $8$ as an upper bound for $C_3$, instead of $6$. Now, this idea of using multisets instead of sets can actually be used more generally to lower $u_k$ and thereby $C_p$ in certain cases. \\

To give just one example, with $p = 2^{20} - 3$, we may take $S = S_1$ to be equal to the set $\{p, 2^1, 2^2, \ldots, 2^{20}\}$, by the proof of Lemma \ref{ssize}. This gives $M_2 = |S| = 21$, so that all coefficients $c_i$ in the second interval are at most $21$. In step $i$ with $v_1 < i \le v_2$, our algorithm as described above writes $c_i$ or $c_i - 1$ as a sum of powers of two, using $S_2 := \{2, 4, 8, 16\}$. Continuing, we can write $c_i$ or $c_i - 1$ in the third interval as a subset sum of $S_3 := \{2, 4 \}$, as $c_i$ is at most $4$ there, while for the fourth interval we use $S_4 := \{2\}$. When we work this out, we obtain $M_3 = 4$, $M_4 = 2$ and $M_5 = 1$ with $u_2 = p+3$, $u_3 = 16$, $u_4 = 4$ and $u_5 = 2$. All in all we deduce that, for this value of $p$, we can take $C_p = 128(p+3)$. This can be optimized, however. \\

Instead of using the sets $S_1$, $S_2$, $S_3$ and $S_4$ to write the coefficients in the first four intervals, we claim that we are better off using the multisets
\begin{align*}
S'_1 &:= \{p, 2, 2, 4, 8, 8, 2^4, 2^5, \ldots, 2^{18}, 2^{19}, 2^{19} \}, \\ 
S'_2 &:= \{2, 2, 4, 8, 8\}, \\
S'_3 &:= \{2, 2\}  \text{ and} \\
S'_4 &:= \{2\}.
\end{align*}

Now, $|S'_1|$ and $|S'_2|$ are larger than $|S_1|$ and $|S_2|$ respectively, which implies that the coefficients in the first three intervals may increase as well. In general, if we write the coefficients in $I_{k-1}$ with `distinct' elements from $S'_{k-1}$ and assume $S'_k \subset S'_{k-1}$, then the coefficients in $I_k$ can be as large as $|S'_{k-1}|$. That is, the upper bound (\ref{maxmk}) changes from $c_i \le M_k$ to $c_i \le |S'_{k-1}|$. On the other hand, the $u_k$ can be redefined to $\max S'_{k-1}$, which improves the previous bound to $C_p = 32p$. In fact, even these improved sets are not optimal, as $S''_1 = \{2, 4, p, p+3 \}$ would lead to an even better constant. \\

It is unclear whether such considerations can lead to significantly lowering the upper bound on $C_p$ in general. In particular, it remains an interesting challenge to show the existence of a constant $c$ such that $C_p < cp$ holds for all odd $p > 1$.

\end{document}